\documentclass{elsart}
\journal{Journal of Algebra}
\usepackage{amsmath}
\usepackage{amsfonts}
\usepackage{amssymb}

\DeclareMathOperator{\rank}{rk}

\DeclareMathOperator{\ord}{ord}
\DeclareMathOperator{\charact}{char}

\newcommand{\Le}{\leqslant}
\newcommand{\Ge}{\geqslant}

\newcommand{\K}{{\mathbf k}}
\newcommand{\ZN}{{\mathbb Z}_{\Ge 0}}

\newcommand{\Z}{{\mathbb Z}}
\newcommand{\C}{{\mathcal C}}
\newcommand{\A}{{\mathcal A}}

\def\B{{\mathcal B}}

\newcommand{\ld}{\mathop{\rm ld}\nolimits}

\newcommand{\rk}{\mathop{\rm rk}\nolimits}

\newcommand{\algrem}{\mathop{\sf algrem}\nolimits}

\newcommand{\sep}{{\;|\;}}

\renewcommand{\u}{{\bf u}}
\renewcommand{\i}{{\bf i}}
\newcommand{\s}{{\bf s}}

\renewcommand{\b}{\ \ \ }

\makeatletter
\def\Ddots{\mathinner{\mkern1mu\raise\p@
\vbox{\kern7\p@\hbox{.}}\mkern2mu
\raise4\p@\hbox{.}\mkern2mu\raise7\p@\hbox{.}\mkern1mu}}
\makeatother

\begin{document}

\begin{frontmatter}
\title{A Bound for Orders in Differential Nullstellensatz\thanksref{thank1}}
\author[OG]{Oleg Golubitsky\thanksref{thank3}}
\ead{Oleg.Golubitsky@gmail.com}
\author[MVK]{Marina Kondratieva}
\ead{kondrmar@rol.ru}
\author[AO]{Alexey Ovchinnikov\thanksref{thank2}}
\ead{aiovchin@math.uic.edu}
\author[AS]{Agnes Szanto\thanksref{thank4}}
\ead{aszanto@ncsu.edu}
\address[OG]{University of Western Ontario\\Department of Computer Science\\London, Ontario, Canada N6A 5B7}
\address[MVK]{Moscow State University\\ Department of Mechanics and Mathematics\\ Leninskie gory, Moscow, Russia, 119991}
\address[AO]{University of Illinois at Chicago\\ Department of Mathematics, Statistics, and Computer Science\\ Chicago, IL 60607-7045, USA}
\address[AS]{North Carolina State University\\Department of Mathematics\\ Raleigh, NC 27695-8205, USA}
\date\today
\thanks[thank1]{The work was partially supported by the Russian Foundation for Basic Research, project no. 05-01-00671.}
\thanks[thank2]{This author was also partially supported by NSF Grant CCR-0096842.}
\thanks[thank3]{This author was also partially supported by NSERC Grant PDF-301108-2004.}
\thanks[thank4]{This author was also partially supported by NSF Grant CCR-0347506}

\begin{keyword}
differential algebra \sep characteristic sets \sep radical differential ideals \sep differential Nullstellensatz
\MSC 12H05 \sep 13N10 \sep 13P10
\end{keyword}

\begin{abstract}
We give the first known bound for orders of differentiations in differential Nullstellensatz for both partial and ordinary algebraic differential equations. This problem was previously addressed in~\cite{Seidenberg} but no complete solution was given. Our result is a complement to the corresponding result in algebraic geometry,
which gives a bound on degrees of polynomial coefficients in
effective Nullstellensatz~\cite{Hermann,Mayr,Brownawell,Kollar,Galligo1,Galligo2,Teresa,Jalonek}.

\medskip
This paper is dedicated to the memory of Eugeny Pankratiev, who was the advisor of the
first three authors at Moscow State University.
\end{abstract}

\end{frontmatter}

\section{Introduction}
Given a system of algebraic partial differential equations $F = 0$, where $F=f_1,\ldots,f_k$, and a differential equation $f=0$, one can effectively
test if $f$ is a differential algebraic consequence
of $F$. In this paper we develop a method that
leads to an effective procedure which finds an
algebraic expression of some power of $f$ in terms of the elements of $F$ and their derivatives  (or shows that such an expression does not exist). This procedure is called {\it effective differential Nullstellensatz}.
A brute-force algorithm solving this problem
consists of two steps:
\begin{enumerate}
\item find an upper bound $h$ on the number of differentiations one needs to apply to $F$ and
\item find an upper bound on the degrees of polynomial coefficients
$g_i$ and a positive integer $k$
\end{enumerate}
such that $f^k$ is a combination of the elements of $F$ together with
the derivatives up to the order $h$ and the coefficients $g_i$.
We solve the first problem in the paper. The second problem was
addressed and solved in \cite{Hermann} and further analyzed and
improved in \cite{Mayr,Brownawell}. A purely algebraic solution was
given in \cite{Kollar}. Most of the references on the subject can
be found in~\cite{Galligo1,Galligo2,Teresa,Jalonek}.

More precisely, our problem is as follows.
We are given a finite set $F$ 
of differential polynomials such that 
a differential polynomial $f$ belongs
to the radical differential ideal generated by $F$ in the
ring of differential polynomials.
Knowing {\bf only} the orders and degrees of the elements of $F$ and the order of $f$, we find a non-negative integer $h$ such that $f$ belongs to the radical of the algebraic ideal generated by $F$ and its derivatives up to the order $h$.

We give a complete solution to this problem using differential
elimination. The problem is non-trivial: the first (unsuccessful) attempt  was made by
Seidenberg \cite{Seidenberg}, where it was conjectured that most
likely such a bound would not be found. Here is where the main
difficulty is coming from.
In order to get the bound  using  a differential elimination algorithm we need to
estimate how many differentiation steps this
algorithm makes. Originally, termination proofs
for such algorithms were based on the Ritt-Noetherianity of the ring of differential polynomials, that is: every increasing chain
of radical differential ideals terminates. And
this result does not say when the sequence terminates. We overcome this problem in
our paper.

The article is organized as follows. We introduce basic notions of differential algebra in Section~\ref{Basics}. Then  we formulate
the main result, Theorem~\ref{MainTheorem}, in Section~\ref{MainResult}. In order to achieve this,
we bound the length
of  increasing sequences of radical differential ideals appearing in our differential elimination in Section~\ref{EliminationAlgorithm} (see Proposition~\ref{algbound}).
For that, in Section~\ref{BoundsForSequences}, we first bound the length of dicksonian
    sequences of tuples of natural numbers with restricted growth
    of the maximal element in these tuples (Lemma~\ref{lexbound}).
Finally, we apply this to obtain the bound for
the differential Nullstellensatz in Theorem~\ref{MainTheor}, from which Theorem~\ref{MainTheorem} follows.
 We conclude by giving in Section~\ref{MichaelsProof} an alternative non-constructive
proof of existence of the bound, based on model theory.

There is some previous work on bounding orders in differential elimination
algorithms. In the ordinary case, we can bound the orders of derivatives of
the output and all intermediate steps of differential elimination
\cite{Bounds}, and this
bound holds for any ranking. Also in the ordinary case, one can give bounds
for quantifier elimination~\cite{Grigoriev} and for the orders and degrees of resolvents
of prime differential ideals of a certain type~\cite{Solerno}. A related bound
for involutive
prolongation, based on the analysis of stability of Spencer sequences,
is obtained in~\cite{Sweeney}.

Note that, unlike the bounds for differential elimination mentioned above, the
bound for the differential Nullstellensatz proposed in this paper holds for the
PDE case. Our bound is also based on the analysis of differential elimination.
But, due to the ranking-independent nature of the differential Nullstellensatz,
we could restrict our analysis to orderly rankings, which allowed us to treat
not only the ordinary case, but the PDE case as well.


\section{Basic differential algebra}\label{Basics}

One can find recent tutorials on the constructive theory of
differential ideals in \cite{Dif,Sit,BouChaptire}. One also refers
to~\cite{Seidenberg,Rit,Kol,Fac,Imp,Bou1,Bou2,Bou3,CF3,pan,Und} for differential
elimination theory. A differential ring is a commutative ring with unity endowed with a set of derivations $\Delta =
\{\partial_1,\ldots,\partial_m\}$, which commute pairwise. The case of
$\Delta = \{\delta\}$ is called {\it ordinary}. 
 Construct the multiplicative monoid
$\Theta = \left\{\partial_1^{k_1}\partial_2^{k_2}\cdots\partial_m^{k_m}\;\big|\; k_i \Ge 0\right\}$
of {\it derivative operators}. Let $Y=\{y_1,\ldots,y_n\}$ be a set whose
elements are called {\it differential indeterminates}.
The elements of the set $\Theta Y=\{\theta y\;|\;\theta\in\Theta,\;y\in Y\}$
are called {\it derivatives}. Derivative operators from $\Theta$ act on
derivatives as
$\theta_1(\theta_2y_i) = (\theta_1\theta_2)y_i$ for all
$\theta_1,\theta_2 \in \Theta$ and $1 \Le i \Le n$.

The ring of {\it differential polynomials} in differential
indeterminates $Y$ over a differential field $\K$ is a ring of commutative
polynomials with coefficients in $\K$ in the infinite set of variables
$\Theta Y$. This ring is denoted by $\K\{y_1,\dots,y_n\}$. We consider the case
of $\charact \K=0$ only. An ideal $I$ in $\K\{y_1,\ldots,y_n\}$ is called {\it
  differential}, if for all $f\in I$ and $\delta\in\Delta$, $\delta
f\in I$.
Let $F \subset \K\{y_1,\dots,y_n\}$ be a set of differential
polynomials. For the differential and radical differential
ideal generated by $F$ in $\K\{y_1,\dots,y_n\}$, we use
notations $[F]$ and $\{F\}$, respectively.

A {\it ranking} 
is a total order $>$ on the set $\Theta Y$
satisfying the following conditions for all $\theta\in\Theta$ and
$u,v\in\Theta Y$:
\begin{enumerate}
\item $\theta u \Ge u,$
\item $u \Ge v \Longrightarrow \theta u \Ge \theta v.$
\end{enumerate}
Let $u$ be a derivative,
that is, $u = \theta y_j$ for
$\theta = \partial_1^{k_1}\delta_2^{k_2}\cdots\partial_m^{k_m} \in \Theta$
 and $1\Le j \Le n$. The {\it order} of $u$ is defined as
$$\ord u=\ord\theta=k_1+\ldots+k_m.$$ If $f$ is a differential
polynomial, $f\not\in\K$, then $\ord f$ denotes the maximal order of
derivatives appearing effectively in $f$.

A ranking $>$ is called
{\it orderly} if $\ord u > \ord v$ implies $u > v$ for all
derivatives $u$ and $v$. 
Let a ranking $>$ be fixed.
The derivative $\theta y_j$ of the highest rank appearing
in a differential polynomial $f \in \K\{y_1,\dots,y_n\} \setminus \K$
is called the {\it leader} of $f$. We denote the leader by $\ld f$ or $\u_f$.
Represent $f$ as a univariate polynomial in $\u_f$:
$$
f = \i_f \u_f^d + a_1 \u_f^{d-1} + \ldots + a_d.
$$
The monomial $\u_f^d$ is called the {\it rank} of $f$ and
is denoted by $\rk f.$
Extend the ranking relation on derivatives  to ranks:
$u_1^{d_1}>u_2^{d_2}$ if either $u_1>u_2$ or
$u_1=u_2$ and $d_1>d_2$.
The polynomial $\i_f$ is called the {\it initial} of $f$.
Apply any derivation $\delta \in \Delta$ to $f$:
$$ \delta f = \frac{\partial f}{\partial \u_f}\delta \u_f + \delta \i_f \u_f^d +
\delta a_1 \u_f^{d-1}+\ldots + \delta a_d.
$$
The leader of $\delta f$ is $\delta \u_f$ and the initial
of $\delta f$ is called the {\it separant} of $f$, denoted $\s_f$.
If $\theta\in\Theta\setminus\{1\}$, then $\theta f$ is called a
{\it proper derivative} of $f$. Note that the initial of any
proper derivative of $f$ is equal to $\s_f$.

We say that a differential polynomial $f$ is {\it partially reduced}
w.r.t. $g$ if no proper derivative of $\u_g$ appears in $f$.
A differential polynomial $f$ is {\it algebraically reduced} w.r.t.
$g$ if $\deg_{\u_g} f<\deg_{\u_g}g$.
A differential polynomial $f$ is {\it reduced} w.r.t. a differential
 polynomial $g$ if $f$ is partially and algebraically reduced
 w.r.t. $g$. Consider any subset
$\A \subset \K\{y_1,\ldots,y_n\}\setminus \K$. We say that
$\A$ is autoreduced (respectively, algebraically autoreduced)
if each element of $\A$ is reduced (respectively, algebraically reduced)
w.r.t. all the others.

Every autoreduced set is finite
\cite[Chapter I, Section 9]{Kol} (but an algebraically autoreduced set
in a ring of differential polynomials
may be infinite). For autoreduced sets we use capital calligraphic letters
$\mathcal{A, B, C,}$ \ldots and
notation $\A = A_1,\ldots,A_p$ to specify the list of the elements
of $\A$ arranged in order of increasing rank.
We denote the sets of initials and separants of
elements of $\A$ by $\i_\A$ and $\s_\A$,
 respectively. Let $H_\A=\i_\A\cup \s_\A$.
For a finite set $S$  of differential polynomials denote 
by $S^\infty$ the multiplicative set containing $1$ and
 generated by $S$. Let $I$ be an ideal in a commutative ring $R$.
The {\it saturated ideal} $I:S^\infty$ is defined as
$\{a \in R\:|\:\exists s \in S^\infty: sa \in I\}$. If $I$ is a
differential ideal then $I:S^\infty$ is also a differential ideal
(see \cite{Kol}).

Let $\A = A_1,\ldots,A_r$ and $\B = B_1,\ldots,B_s$
 be (algebraically) autoreduced sets. We say that $\A$ has lower rank than
$\B$ if
\begin{itemize}
\item
there exists $k \Le \min(r, s)$ such that $\rank A_i$ =
$\rank B_i$ for $1 \Le i < k,$ and $\rank A_k < \rank B_k$,
\item or if $r > s$ and $\rank A_i = \rank B_i$ for $1 \Le i \Le s$.
\end{itemize}
We say that
 $\rank\A = \rank\B$ if $r=s$ and
$\rank A_i = \rank B_i$ for $1 \Le i \Le r$.
Let $v$ be a derivative in $\K\{y_1,\ldots,y_n\}.$ Denote by $\A_v$ the set of the elements of $\A$ and their derivatives that have a leader ranking strictly lower than $v$.
A set $\A$ is called {\it coherent} if whenever $A, B \in \A$ are such that $\u_A$ and $\u_B$ have a common derivative: $v = \psi \u_A = \phi\u_B$, then $\s_B\psi A - \s_A\phi B \in (\A_v):H_\A^\infty$.

\section{Main result}\label{MainResult}

For a finite set of differential polynomials $F \subset
\K\{y_1,\ldots,y_n\}$ let $D(F)$ be the maximal total degree of a polynomial in $F$.
For each $i$, $1 \Le i\Le n$, let
$$
h_i(F) = \ord_{y_i}(F),\quad
H(F) = \max_{1\Le i\Le n}h_i(F).
$$
For $h \in \Z_{\Ge 0}$ let $F^{(\Le h)}$ denote the set of derivatives of the elements of $F$ of order less than or equal to $h$.
The {\it Ackermann function} appearing in our main result is defined as follows \cite[Section 2.5.5]{Cut}:
\begin{align*}
   A(0,n)&=n+1\\
   A(m+1,0)&=A(m,1)\\
   A(m+1,n+1)&=A(m,A(m+1,n)).
  \end{align*}
\begin{thm}\label{MainTheorem} Let $F \subset \K\{y_1,\ldots,y_n\}$ be a finite set, $0 \ne f \in \{F\}$ and
let $t(F,f)$ be the minimal non-negative integer such that
$f \in\sqrt{\left(F^{(\Le t(F,f))}\right)}$. 
Then
$$
t(F,f) \Le A(m+8,\max(n,H(F\cup f),D(F\cup f))).
$$
\end{thm}
\begin{pf}
This result will be proved step-by-step in the following sections as described in the introduction and  finally established in Theorem~\ref{MainTheor}. 
\end{pf}

\begin{rem}It is our own choice here to bound $t(F,f)$ using solely the maximal
orders and degrees of $F$ and $f$. One might come up with another
bound using more information of $F$ and $f$. But we emphasise that the bound on orders {\bf must} depend on the degrees, number of differential indeterminates,
and number of basic differentiations as the following examples
show.
\end{rem}

\begin{exmp} Let $f = 1$ and $F = \left\{y'-1,y^k\right\}$ in $\K\{y\}$, the ordinary case. In order
to express $1$ in terms of the elements of $F$, one has to differentiate $y^k$ $k$
times.
\end{exmp}

In the linear and non-linear cases, consider the following examples showing that the bound must depend on the number of variables and derivations.

\begin{exmp} Let $f = 1$ and $F = \left\{y_1',y_1-y_2',\ldots,y_{n-1}-y_n',y_n-a\right\}$ in the
ordinary differential ring $\K\{y_1,\ldots,y_n\}$,
where $a \in \K$ is such that $a^{(n)} = 1$. We have to differentiate the first $n-1$ generators $n$ times to get $y_n^{(n)}$ into the corresponding algebraic ideal. Hence, $t(F,f) = n$.
\end{exmp}

\begin{exmp} Let $f =1$ and $F = \{y_1^2,y_1-y_2^2,\ldots,y_{n-1}-y_n^2, 1-y_n'\} \subset\K\{y_1,\ldots,y_n\}$, again in the ordinary case. 
One can show that
\begin{align*}
F\subset&\left(y_1,y_2,\dots,y_n,1-y_n'\right)=I_0,\\
F^{(\Le 1)}\subset&\left(I_0,y_1',y_2',\dots,y_{n-1}',y_n''\right)=I_1,\\
F^{(\Le 2)}\subset&\left(I_1,y_1'',\dots,y_{n-2}'',y_{n-1}''-2,y_n^{(3)}\right)=I_2,\\
F^{(\Le 3)}\subset&\left(I_2,y_1''',\dots,y_{n-2}''',y_{n-1}''',y_n^{(4)}\right)=I_3,\\
F^{(\Le 4)}\subset&\left(I_3,y_1^{(4)},\dots,y_{n-2}^{(4)}-2^2\binom{4}{2},y_{n-1}^{(4)},y_n^{(5)}\right)=I_4,\\
&\ldots\\
F^{\Le(2^{n-1})}\subset&\left(I_{2^{n-1}-1},y_1^{(2^{n-1})}-\prod_{k=1}^{n-1}\binom{2^k}{2^{k-1}}^{2^{n-k-1}},y_2^{(2^{n-1})},\ldots,y_n^{\left(2^{n-1}+1\right)}\right)=I_{2^{n-1}},\\
&\ldots\\
F^{\Le(2^n-1)}\subset&\left(I_{2^n-2},y_1^{(2^n-1)},y_2^{(2^n-1)},\ldots,y_n^{(2^n)}\right).
\end{align*}
Therefore, $1 \notin\left(F^{\Le (2^n-1)}\right)$.
Thus, $t(F,f) = 2^n$, because modulo $1-y_n'$ we have 
$$
\left(y_n^{2^n}\right)^{(2^n)} = 2^n\left(\left(y_n^{2^n-1}\right)y_n'\right)^{(2^n-1)} \equiv 2^n\left(y_n^{2^n-1}\right)^{(2^n-1)}\equiv\ldots\equiv(2^n)!(y_ny_n')'\equiv 2^n!y_n'^2\equiv1.
$$
\end{exmp}

\begin{exmp} If we replace $F$ in the previous example  by 
$$G = \left\{u_{x_1}^2,u_{x_1}-u_{x_2}^2,\ldots,u_{x_{m-1}}-u_{x_m}^2,1-u_{x_m^2}\right\} \subset\K\{u\}$$ with partial derivatives $\partial_{x_1},\ldots,\partial_{x_m}$, we obtain
an example which shows that the bound on orders must depend on the number $m$
of derivations. Again, the generators will have to be differentiated $2^m$ times
to express $1$. 
\end{exmp}

\section{Bounds on lengths of sequences}\label{BoundsForSequences}
The results of this section with be further used in Section~\ref{EliminationAlgorithm} to bound lengths of decreasing sequences of autoreduced sets
appearing in the differential elimination algorithm that we use. In this section the letters $m$ and
$n$ will {\it not} mean the number of derivations and differential indeterminates,
respectively.

We begin by bounding the length of certain sequences of non-negative $n$-tuples.
 Call a sequence $t_1,t_2,\ldots,t_k$ of $n$-tuples 
{\it dicksonian}, if for all $1\Le i<j\Le k$, there does not exist a 
non-negative $n$-tuple $t$ such that $t_i+t=t_j$. 
For example, any lexicographically decreasing sequence is dicksonian.
By Dickson's Lemma, every dicksonian sequence is finite. Our goal 
is to obtain an explicit upper bound for the length of a dicksonian
sequence, whose elements do not grow faster than a given function,
in terms of this function, the first element, and the size $n$ of 
the tuples.
Let
$$(a_1^1,\ldots,a_n^1),(a_1^2,\ldots,a_n^2),\ldots
,(a_1^k,\ldots,a_n^k)
$$
be a dicksonian sequence of $n$-tuples of non-negative integers such that
\begin{align}\label{ineq1}
\max\left(a_1^j,\ldots,a_n^j\right)\Le f(j)
\end{align}
for all $j$, $1\Le j\Le k$, where
$$
f: \Z_{\Ge 0} \to \Z_{\Ge 0}
$$
is a fixed function.
We say that the {\it growth of this sequence is bounded by the function $f$}. 

The following proposition closely resembles a particular case of our problem, 
namely that of $f(i)=m+i-1$. However, in Proposition~\ref{p1} the maximal coordinate 
must increase by 1 at each step, whereas in our case it is allowed
to decrease or remain the same. We will reduce the case of 
a dicksonian sequence with the growth bounded by a function $f$
from a certain large class of functions that ``do not grow too fast'',
to the one treated in Proposition~\ref{p1}. 

\begin{prop}\cite[Proposition 1]{SocAck} \label{p1}
   Let $t_1,t_2,\ldots,t_k$ be a dicksonian sequence of $n$-tuples, such 
   that the maximal coordinate of $t_i$ equals $m+i-1$, for all $1\Le i\Le k$.
   Then the maximal coordinate in the last tuple, $t_k$, does not exceed $A(n,m-1)-1$,
   and there exists such a dicksonian sequence for which this bound is reached.
\end{prop}

Note that in Proposition~\ref{p1} we have: $m$ is the maximal coordinate of $t_1$ and
the length $k$ of the sequence is bounded by $A(n,m-1)-m$.
The general case (of any function $f$, not necessarily 
from our class) has also been studied in \cite{SocLen} using a different approach. 
It is shown that the maximal possible length is primitive recursive in $f$ and 
recursive, but not primitive recursive (if $f$ increases at least linearly), 
in $n$. Sequences yielding the maximal possible length are constructed. Moreover, 
if $f$ is linear, an explicit expression for the maximal length is given in terms 
of a generalized Ackermann function. Our statement was motivated by the need to 
obtain an explicit expression for the bound for a wider class of growth 
functions.

Let $L_{f,n}$ denote the maximal length of a dicksonian sequence of $n$-tuples, whose growth
is bounded by $f$.
For an increasing function $f:\ZN\to\ZN$, let $\lceil f^{-1}(x)\rceil$ be the least number $k$
such that $f(k)\Ge x$.

\begin{lem}\label{lexbound} Let $f:\ZN \to \ZN$ be an increasing function,
and $d\in\ZN$ be a number such that $f(i+1)-f(i)\Le A(d,f(i)-1)$ for all 
$i>0$. Then 
\begin{align}\label{ineq2}
L_{f,n} < \big\lceil f^{-1}\big(A(n+d,f(1)-1)\big)\big\rceil
\end{align}
and the maximal entry of the last $n$-tuple does not exceed $A(n+d,f(1)-1)$. 
\end{lem}
\begin{pf}
Consider a disckonian sequence 
\begin{align}\label{seq1}
(a_1^1,\ldots,a_n^1),\;
(a_1^2,\ldots,a_n^2),\ldots,\;
(a_1^k,\ldots,a_n^k),
\end{align}
whose growth is bounded by $f$.
Construct from \eqref{seq1} a new sequence satisfying the conditions of 
Proposition~\ref{p1}. 
Append to the first tuple $d$ new coordinates, each equal to $f(1)$,
obtaining the following $(n+d)$-tuple:
$$(a_1^1,\ldots,a_n^1,f(1),\ldots,f(1)).$$
Then add $f(2)-f(1)-1$ new $(n+d)$-tuples. 
The first $n$ coordinates of these tuples are
$(a_1^1,\ldots,a_n^1)$. The last $d$ coordinates form a dicksonian
sequence of $d$-tuples, starting with $(f(1),\ldots,f(1))$, with the 
maximum coordinate growing exactly by 1 at each step. 
From Proposition~\ref{p1} and condition
$f(2)-f(1)\Le A(d,f(1)-1)$, such sequence exists. The last tuple will 
have the maximum coordinate equal to $f(2)-1$.
Next, add the tuple
$$(a_1^2,\ldots,a_n^2,f(2),\ldots,f(2)).$$
Since the growth of \eqref{seq1} is bounded by $f$, the maximal
coordinate in this tuple equals $f(2)$.
Continue by adding $f(3)-f(2)-1$ new $(n+d)$-tuples, whose
first $n$ coordinates are $(a_1^2,\ldots,a_n^2)$ and last
$d$ coordinates form a dicksonian sequence growing by 1 at
each step.
Finally, when the tuple
$$(a_1^k,\ldots,a_n^k,f(k),\ldots,f(k))$$
is reached, stop. 
We obtain a sequence of $(n+d)$-tuples in which the 
maximal coordinate grows by 1 at each step. We will show
that this sequence is dicksonian.
Suppose that it is not. 
Let $t_j$, $t_l$, $j<l$, be two $(n+d)$-tuples from this sequence, 
for which there exists a tuple $t$ such that $t_l=t_j+t$. 
Let $t^{I}$, $t^{II}$ denote the first $n$ 
coordinates and the last $d$ coordinates of an $(n+d)$-tuple $t$, 
respectively. Then we have $t_l^{I}=t_j^I+t^I$ and 
$t_l^{II}=t_j^{II}+t^{II}$.
If $t_j$ and $t_l$ have been added after {\it the same} tuple 
of the form $$p_i=(a_1^i,\ldots,a_n^i,f(i),\ldots,f(i)),$$
or if $t_j$ coincides with such a tuple $p_i$ and $t_l$ has been 
added after $p_i$, 
the equality $t_l^{II}=t_j^{II}+t^{II}$ contradicts the fact that the 
last $d$ coordinates of the tuples between $p_i$ and $p_{i+1}$,
including $p_i$ and excluding $p_{i+1}$, form a dicksonian sequence.
If $t_j$ and $t_l$ have been added after {\it different} tuples
$p_i$ and $p_{i'}$, the equality
$t_l^{I}=t_j^I+t^I$  contradicts the fact that sequence
\eqref{seq1} is dicksonian.
Therefore, our assumption was false and the constructed sequence
is dicksonian. 

By Proposition~\ref{p1}, the maximum coordinate
of its last element does not exceed $A(n+d,m-1)-1$. Since the
maximum coordinate in the first element is $f(1)$ and grows by 1 at each
step, the number of elements in the constructed sequence does not 
exceed $A(n+d,f(1)-1)-f(1)$. 
On the other hand, the number of elements in the constructed
sequence is:
$$f(2)-f(1)+f(3)-f(2)+\ldots+f(k)-f(k-1)+1=f(k)-f(1)+1.$$
Therefore,
$$f(k)-f(1)+1\Le A(n+d,f(1)-1)-f(1),$$
that is, 
$$f(k)< A(n+d,f(1)-1),$$
and 
$$k<\big\lceil f^{-1}\big(A(n+d,f(1)-1)\big)\big\rceil.$$
\end{pf}

\section{Differential elimination algorithm}\label{EliminationAlgorithm}
Using the result of the previous section, we obtain an upper bound 
for the length of sequences of autoreduced sets of decreasing rank 
produced by a differential elimination algorithm. The idea is to 
put in correspondence with such a sequence a dicksonian sequence 
of tuples, whose growth is bounded by a function derived from the 
algorithm.

We fix an {\bf orderly} ranking.
Algorithm~\ref{MainRC} computes a characteristic decomposition of a radical
differential ideal given by a set of generators. It is designed in such a way that
allows us to control the orders and degrees of differential polynomials occurring
in all intermediate steps together with a {\bf bound} on the number of iterations
of this algorithm.
Also, in the algorithm the procedure $\algrem$ computes an algebraic pseudo-remainder of a polynomial with respect to an algebraic triangular set (a set is called {\it triangular} if the leaders of its
elements are distinct
).
Algorithm {\sf MinimalTriangularSubset} inputs a finite set of differential polynomials and outputs one of its least rank triangular subsets.
Algorithm {\sf CharSet} inputs a finite set of differential polynomials and outputs one of its characteristic sets, that is, an autoreduced subset of the least rank. Denote by $\Delta(\C)$ the set of ``differential S-polynomials'' of $\C$ defined in~\cite[Definition 4.2]{Dif}.
From now on, $m$ and $n$ again denote the numbers of derivations and differential indeterminates, respectively.

\begin{alg}{\sf RGBound}\label{MainRC}$(F_1)$\\
\begin{tabular}{l}
{\sc Input:} a set 
$F_1 \subset \K\{y_1,\ldots,y_n\}$
with 
derivations $\{\partial_1,\ldots,\partial_m\}$ \\
{\sc Output:} A finite set T of  triangular sets such that $\{F_1\}=\bigcap\limits_{\C\in T}[\C]:H_{\C}^\infty$;\\ 
\hphantom{\sc Output:} 
if $1 \notin [\C]:H_\C^\infty$ then $\C$ is coherent and autoreduced\\
\hphantom{\sc Output:} otherwise $1 \in (\C):H_\C^\infty$.\\
\b $T:=\{\varnothing\};\ U:=\{(F_1, \varnothing)\}$\\
\b {\bf while $U\neq\varnothing$ do}\\
\b \b Take and remove any $(F,\C)\in U$\\ 
\b \b $f:=$ an element of $F$ reduced w.r.t. $\C$ of the least rank\\
\b \b  {\bf if} $\s_f\notin \K$ {\bf then} $U:=U\cup (F\cup \s_f,\C)$ {\bf end if}\\
\b \b  {\bf if} $\i_f\notin\K$ {\bf then} $U:= U\cup (F\cup \i_f,\C)$ {\bf end if}\\ 
\b \b $D:=\{C\in\C\;|\;\ld C=\theta\ld f\;{\rm for}\;{\rm some}\;\theta\in\Theta\}$\\
\b \b $\bar\C:=\C\setminus D\cup\{f\}$;\ \ $G:=F\cup \Delta\left(\bar\C\right)\cup D\setminus\{f\}$\\
\b \b $b:=\max\limits_{g\in G} \ord g$\\
\b \b $\B:=\:${\sf MinimalTriangularSubset}$\left(\left\{\theta C\;|\;C\in\bar\C,\;\ord\theta C\Le b\right\}\right)$\\
\b \b $\bar\B:=\{\algrem(h,\B\setminus\{h\})\;|\;h\in\B\}$\\
\b \b {\bf if} $\rk\bar\B\ne\rk\B$ {\bf then} $T:=T\cup\left\{\B\right\}$; {\bf continue}; {\bf  end if}\\
\b \b $R:=\{\algrem(g,\B)\;|\;g\in G\}\setminus\{0\}$;\ \ $\C:=\:${\sf CharSet}$\left(\bar\B\right)$\\
\b \b {\bf if} $R=\emptyset$ {\bf then} $T:=T\cup\{\C\}$ {\bf else} $U:=U\cup (R\cup F,\C)$ {\bf end if}\\
\b {\bf end while} \\
{\bf return} $T$
\end{tabular}
\end{alg}
We get the following bounds for the growth of the maximal degrees of the polynomials computed at the $i$-th step of Algorithm~\ref{MainRC}. 

\begin{prop}\label{algbound}
Fix $(F_i,\C_i)\neq\varnothing\in U$ and let $(F_{i+1},\C_{i+1})\neq\varnothing$  be any of the
elements obtained from $(F_i,\C_i)$ after one iteration of the {\bf while}-loop.
 We then have
$$
D(F_{i+1}\cup \C_{i+1}) \Le (4D(F_i\cup \C_i))^{\binom{2H(F_i\cup \C_i)+m}{m}+1}.
$$
\end{prop}
\begin{pf}
Consider the  iteration of the loop.
In the first five lines of the loop the degrees do not change as
adding an initial or a separant does not cause an increase in the degrees.
So, the first place where the degrees may change is
the computation of $ \Delta\left(\bar\C\right)$, that is, computing cross differentiations.
This at most doubles the degrees. 
After that, the only places where the degrees
of polynomials may increase are calls to $\algrem(g,\B)$ for $g\in G$
or $\algrem\left(h,\B\setminus\{h\}\right)$ for $h\in \B$.
In both cases it is a sequence of at most $|\B|$ algebraic pseudodivision. We will prove the bound for the reduction of a fixed $g\in G$ modulo $\B$, and the other case follows similarly. 

Assume that $|\B|=N$, and let $\B=\{B_1, \ldots,B_N\}$ be ordered such that 
$$
\ld(B_1)>\ld(B_2)>\cdots>\ld(B_N).
$$
Let $g^{(0)}:=g$ and $g^{(t)}:=\algrem\left(g,\{B_1, \ldots B_t\}\right)$ for $t>0$. 
Denote the maximal total
degree of $g^{(t)}\cup\B$ by $\delta(t)$, $t\Ge0$.
Note that $\delta(0)\Le 2D(F_i\cup\C_i)$. Then
$g^{(t+1)}$ is obtained by the pseudo-division of $g^{(t)}$ with respect to the polynomial $B_{t+1}$. Thus,
$$
g^{(t+1)}=\i_{B_{t+1}}^\epsilon g^{(t)} - qB_{t+1},
$$
where $\epsilon$ is a sufficiently large exponent specified below, $q$ is the pseudo-quotient, and the degree of $g^{(t+1)}$ in $\ld(B_{t+1})$ is smaller that the same degree of $B_{t+1}$. 
The exponent $\epsilon$ is bounded by 
$$
\deg_{\ld(B_{t+1})}\left(g^{(t)}\right)-\deg_{\ld(B_{t+1})}(B_{t+1})+1\Le \delta(t).
$$
Therefore, the total degree of $g^{(t+1)}$ is bounded by
$$
\delta(t+1)\Le \delta(t)\delta(0)+\delta(t)+\delta(0)\Le\begin{cases}
 2\delta(0)\delta(t), & t\Ge 1\text{ or } \delta(0)\Ge 2; \\
3,&t=0\text{ and } \delta(0)=1.
 \end{cases}
$$
This implies that 
$\delta(N)\Le (2\delta(0))^{N+1}.$
Using that $\delta(0)\Le 2D(F_i\cup\C_i)$, $N\Le {{b_i+m}\choose{b_i}}$, and $b_i\Le 2H(F_i\cup\C_i)$, we get the claim (the numbers $b_i$ are defined in Algorithm~\ref{MainRC}).
\end{pf}

\subsection{Differential bounds for splitting}\label{splittingsec}
Consider now the splitting part of differential elimination.
Algorithm~\ref{MainRC} removes an element   $(F,\C)$ from $U$ and
within one iteration of the {\bf while}-loop it converts this element into one, two, or
three elements.
Moreover, if the set $\C$ does not change, the orders and degrees of the elements of $F$ do not increase after the conversion.
We call such an iteration {\it incomplete}.
Denote the maximal 
number of all (complete and incomplete) iterations 
of the {\bf while}-loop of Algorithm~\ref{MainRC} 
 by $L(F)$.

\begin{prop}\label{MRCCorrect} Algorithm~\ref{MainRC} is correct and terminates.
Moreover,  
$$
L(F) \Le \log_2(A(m+7,Q(F)-1)),
$$
where
\begin{align}\label{QF}
Q(F)=\max\left(9,n,2^{9H(F)},D(F)\right).
\end{align}
\end{prop}

\begin{pf}To demonstrate {\bf correctness} we will show that
the {\bf while}-loop of Algorithm~\ref{MainRC} has the following invariant
\begin{align}\label{Inv}\{F_1\}=\bigcap_{(F,\:\C)\in
    U}\{F\cup\C\}:H_\C^\infty\bigcap_{\A\in T}[\A]:H_\A^\infty.
    \end{align}
Indeed, since $F_1\subset F$ for all $(F,\C)\in U$ and for any $\A\in T$ 
every element of $F_1$ is reducible to zero with respect to  $\A$, we have
the inclusion ``$\subset$''.
We will show the opposite inclusion by induction on the number of iterations
of the loop (not assuming that it is finite).  
Invariant~\eqref{Inv} holds at the beginning of the first iteration ($T=\varnothing, U=(F_1,\varnothing)$).
Let a finite number of iterations of the loop be executed preserving the invariant. Suppose
that at the next step we remove an element   $(F,\C)$ from $U$.
Let $\bar \C$ be either the set added to $T$ (in this case we let  $\bar F=\emptyset$), or the second element of the pair $\left(\bar F, \bar \C\right)$ returned back to $U$.
Note that in both cases we have 
\begin{align}\label{FCinclude}
\left\{\bar F,\bar\C\right\} \subset \{F,\C\}
\end{align} by construction. 
We will show the inclusion 
$$\left\{\bar F,\bar \C\right\}:H_{\bar\C}^\infty\cap\{F,\i_f,\C\}:H_\C^\infty\cap\{F,\s_f,\C\}:H_\C^\infty\subset \{F,\C\}:H_\C^\infty,$$ that together with the inductive hypothesis shows the result.
Indeed, 
applying~\cite[Proposition 6.6]{Dif}, since the polynomial $f$ chosen in the loop is reduced with respect to $\C$, we have
$$\{F,\C\}:(H_\C\cup H_f)^\infty\cap
  \{F\cup \i_f,\C\}:H_\C^\infty\cap
  \{F\cup \s_f,\C\}:H_\C^\infty=\{F,\C\}:H_\C^\infty.$$
It remains to note that 
$$\left\{\bar F,\bar\C\right\}:H_{\bar\C}^\infty\subset
\{F,\C\}:
\{H_\C
\cup H_f\}^\infty.$$
Indeed,  according to Algorithm~\ref{MainRC}, every $g\in \bar\C$
comes from some $C\in \C$ as a remainder with respect to the  
triangular set $\B\subset[\C]$. Applying~\eqref{FCinclude} together with \cite[Lemma 5]{Bounds} and 
\cite[Lemma 6.9]{Dif} to $K=\s_g$, $H=H_\C$, we obtain that
\begin{align*}
\left\{\bar F,\bar\C\right\}:H_{\bar\C}
&\subset \left\{\bar F,\bar\C\right\}:(H_{\bar \C}\cup  H_\C)^\infty \subset\\
&\subset \{ F,\C\}:(H_{\bar \C}\cup  H_\C)^\infty =\{ F,\C\}:(H_\C\cup H_f)^\infty.
\end{align*}
It remains to note that in the case of $\rk\bar\B \ne \rk\B$ we have $1\in(\B):H_{\B}^\infty$ by~\cite[Lemma 5]{Bounds}.

{\bf Termination} of the algorithm and the {\bf bound} for $L(F)$ will be proved as follows.
To each element of $U$ we associate an $(m+4)$-tuple in such a way that the sequence
of these vectors for each element of  $U$ is dicksonian.
Note that by our definition  $L(F)$ is the maximal length of such dicksonian sequence. 
At the beginning, the set $U$ consists of one element only. We associate the vector 
$$\tau_0=(0,\ldots,0,n,n,D(F))$$ to it.
Let an element  $(F,\C)$ with vector $\tau$ be under processing of the loop and let $f$ 
be an element of $F$ reduced with respect to  $\C$ and of minimal rank with this property. 

If an {\it incomplete} iteration occurs, we add an element of the form $(F\cup\{g\},\C)$ to $U$. The degree
of  $g$ is less than the one of $f$ (this is either the initial or separant of $f$). We then transform the vector $\tau$ by replacing the last its component by $D(g)$. 
If the iteration is {\it complete} then $U$ is concatenated with a converted element $\left(\bar F,\bar\C\right)$. In this case we change the components of $\tau$ in the following way.  
  Let $\rk f=(\theta y_j)^d$ be the rank of $f$,  
  where $\theta=\partial_1^{i_1}\ldots\partial_m^{i_m}$.
  Consider the $(m+4)$-tuple
$$\bar\tau=(i_1,\ldots,i_m,d,j,n-j,\deg(g)),$$
where $g$ is the element of  $\bar F$ reduced with respect to $\bar \C$ 
and of minimal rank. Proceeding in the described way from $\tau$ to $\bar\tau$, let $\tau_k$ denote the $(m+4)$-tuple corresponding to the $k$-th iteration of the {\bf while}-loop.
  Since $ f$ is reduced w.r.t. $\C$,
  the sequence $\tau_0,\tau_1,\tau_2,\ldots$ is dicksonian. Therefore, the number of
  iterations for each element is bounded.
This proves {\bf termination} of
  Algorithm~\ref{MainRC}. Note that if we remove $n-j$ from $\bar\tau$, the sequence might not be dicksonian. Indeed, let  $F=\{y_1,y_2\}\subset \K\{y_1,y_2,y_3\}$ with $y_1<y_2<y_3$. We then would have $\tau_1=(0,1,1,1)$ and $\tau_2=(0,1,2,1)$.
  
Let $H_1=\max(H(F_1),m)$, $H_{k+1}=2H_k$, $D_1=D(F_1)$, and 
  $$D_{k+1}=(4D_k)^{\binom{2H_k+m}{m}+1}.$$
  By Proposition~\ref{algbound}, the maximal coordinate
  of the $(m+4)$-tuple $\tau_k$ does not exceed $\max(H_k,D_k,n)$.
    Let 
$$
u_1=\max\left(n,9,2^{9H(F)},D(F)\right),\;\;\;
        u_{k+1}=2^{\sqrt[3]{u_k}(2+\log_2u_k)}.
$$
  Then the maximal coordinate of $\tau_k$ does not exceed $u_k$
  for all $k\Ge 1$. Indeed, 
  we will prove by induction that $H_k \Le \frac {1}{9}\log_2 u_k$ and 
  $D_k\Le u_k$.
  For $k=1$ these inequalities hold by definition of $u_1$.
  Assuming that they hold for $H_k$, $D_k$, and $u_k$, prove them for $H_{k+1}$,
  $D_{k+1}$, $u_{k+1}$. Since $H_{k+1}=2H_k$, we have 
  \begin{align*}
2^{9H_{k+1}}=2^{9\cdot 2H_k}=\left(2^{9H_k}\right)^2\Le 
    u_k^2 = 2^{2\log_2u_k} \Le 2^{\sqrt[3]{u_k}\log_2u_k} < u_{k+1},
    \end{align*}
  because $u_k\Ge 9$. 
  Next,
  \begin{align*}
\log_2 D_{k+1}&=\left(\binom{2H_k+m}{m}+1\right)(2+\log_2 D_k)\Le {2^{2H_k+m}}(2+\log_2 u_k)\Le\\
  &\Le 2^{3H_k}(2+\log_2 u_k)\Le \sqrt[3]{u_k}(2+\log_2 u_k)=\log_2 u_{k+1}.
  \end{align*}
  Here we used the fact that $H_k\Ge m$, as well as the inequality
  $H_k\Le \frac{1}{9}\log_2 u_k$ proven above.
    Now observe that for $x\Ge 9$, we have
  $$2^{\sqrt[3]x(2+\log_2 x)}\Le 2^{x+2}-3=A(3,x-1).$$ 
  Therefore, the sequence of $(m+4)$-tuples $\tau_0,\tau_1,\tau_2,\ldots$ 
  satisfies the conditions of Lemma~\ref{lexbound} with $d=3$. And, according to this lemma,
  the length of this sequence does not exceed
  $$\big\lceil f^{-1}\big(A(m+7,f(1)-1)\big)\big\rceil,$$
  where $f(k)=u_k \Ge 2^k$. We can now plug in $f(1)=u_1$, and replace $f^{-1}$ with $\log_2$.
\end{pf}
\begin{cor}\label{maxcor} The maximal orders
and degrees of polynomials computed by Algorithm~\ref{MainRC} do not exceed
\begin{align}\label{boundformula}
 A(m+7,Q(F)-1).
\end{align}
\end{cor}
\begin{pf} 
Follows directly from Lemma~\ref{lexbound} and Proposition~\ref{MRCCorrect}.
\end{pf}

\subsection{Lifting the final bound from splitting}\label{LiftingBound}
Note that for $f \in \K\{y_1,\ldots,y_n\}$ and $F \subset \K\{y_1,\ldots,y_n\}$ we have
$$
1 \in \{F\}:f \iff f\in \{F\}.
$$
\begin{lem}\label{lemmaf} Let $A_1,\ldots,A_p = \A \subset [F]$ be a coherent autoreduced set that reduces all elements of $F$ to zero. Then for all $f \in \{F\}$
$$
1 \in \left(\A^{(\Le q)}\right):\left(H_\A^\infty\cup f\right),
$$
where $q := \max\left(0,\ord f - \min\limits_{g\in \A}\ord g\right)$.
\end{lem}
\begin{pf} By our assumption, $1 \in \{F\}:f$.
Moreover, since $[F] \subset [\A]:H_\A^\infty$, it follows from \cite[Section 5.2]{Fac} that
\begin{align}\label{Decomposition}
\{F\}:f = [\A]:(H_\A^\infty\cup f)\cap\bigcap_{i=1}^p\{F,\i_i\}:f\cap\bigcap_{i=1}^p\{F,\s_i\}:f,
\end{align}
where $\i_i$ and $\s_i$ are the initial and separant of $A_i$, respectively.
Hence, $1 \in [\A]:\left(H_\A^\infty\cup f\right)$, that is, $f \in [\A]:H_\A^\infty$. Let $g$ be a partial pseudo-remainder of $f$
with respect to $\A$. Then by the Rosenfeld lemma \cite[Lemma 5, III.8]{Kol}, we have $g \in (\A):H_\A^\infty$. Since the ranking is orderly, there exsits
$h \in H_\A^\infty$ such that $h\cdot f - g \in \left(\A^{(\Le q)}\right)$. Indeed, at each step of the partial pseudo-division the order of the resulting differential polynomial is less than or equal to the one of the previous polynomial. And $q$ represents the maximal number of times one
possibly needs to differentiate elements of $\A$ to perform one
step of the partial psuedo-reduction. 
\end{pf}

\begin{lem}\label{degreelem} Let $F \subset \K\{y_1,\ldots,y_n\}$ and $f,f_1,\ldots,f_k \in \K\{y_1,\ldots,y_n\}$. Suppose that for some $d\in\Z_{\Ge 1}$ we have $(f_1\cdot\ldots\cdot f_k)^d \in (F):f$. Then 
$$
\theta_1f_1\cdot\ldots\cdot\theta_k f_k  \in \sqrt{\left({F^{\left(\Le4^{(k+1)H+1}d\right)}}\right)}:f,
$$
where $\theta_i \in \Theta$ with $\ord\theta_i \Le H$ for all $i$, $1 \Le i \Le k$.
\end{lem}
\begin{pf}It follows from the proof of \cite[Lemma 1.7]{Kaplansky} that if $a^d \in (F)$ then
$(\partial_i a)^{2d-1}\in \left(F^{(\Le d)}\right)$ for any $a \in \K\{y_1,\ldots,y_n\}$ and $\partial_i \in \Delta$.
Therefore, 
$$
((\partial_if_1)\cdot f_2\cdot\ldots\cdot f_k+\ldots+f_1\cdot\ldots\cdot f_{k-1}\cdot(\partial_if_k))^{2d-1}\in \left(F^{(\Le d)}\right):f^\infty.
$$
Multiplying by $((\partial_if_1)\cdot f_2\cdot\ldots\cdot f_k)^{2d-1}$, we obtain
$$
((\partial_if_1)\cdot f_2\cdot\ldots\cdot f_k)^{2(2d-1)}\in \left(F^{(\Le d)}\right):f^\infty.
$$
To make the computation simpler, we have
$$
((\partial_if_1)\cdot f_2\cdot\ldots\cdot f_k)^{4d}\in \left(F^{(\Le d)}\right):f^\infty.
$$
By induction, we conclude that
$$
((\theta_1 f_1)\cdot f_2\cdot\ldots\cdot f_k)^{4^Hd}\in \left(F^{\left(\Le d\left(1+4+\ldots+4^H\right)\right)}\right):f^\infty \subset  \left(F^{\left(\Le 4^{H+1}d\right)}\right):f^\infty.
$$
Similarly, we obtain
\begin{align*}
((\theta_1 f_1)\cdot (\theta_2f_2)\cdot\ldots\cdot f_k)^{4^{2H}d}\in  &\left(F^{\left(\Le4^{H+1}d+4^{H+1}4^Hd\right)}\right):f^\infty =\\
& = \left(F^{\left(\Le4^{H+1}d(1+4^H)\right)}\right):f^\infty.
\end{align*}
Finally, by induction we get
\begin{align*}
((\theta_1 f_1)\cdot (\theta_2f_2)\cdot\ldots\cdot(\theta_kf_k))^{4^{kH}d}\in & \left(F^{\left(\Le4^{H+1}d\left(1+4^H+\ldots+4^{(k-1)H}\right)\right)}\right):f^\infty
\subset\\
&\subset\left(F^{\left(\Le4^{(k+1)H+1}d\right)}\right):f^\infty ,
\end{align*}
which finishes the proof.
\end{pf}

\begin{lem}\label{NullBound}
For a finite subset $F \subset \K\{y_1,\ldots,y_n\}$ we have:
\begin{align}\label{MainInequality}
t(F,f) \Le \ord f + H(F)\cdot 2^{L(F)} +4^{\left(n\cdot2^{H(F)\cdot2^{L(F)}+1}+1\right)\cdot t(G,f)+1}d,
\end{align}
where 
\begin{align}\label{dformula}
d := \max\left(D(f),A\left(m+7,Q(F)-1\right)\right)^{n\cdot 2^{H(F)\cdot2^{L(F)}+m+\ord f}},
\end{align}
$G \subset \K\{y_1,\ldots,y_n\}$ is such that $L(G) \Le L(F) - 1$ and $H(G) \Le H(F)\cdot 2^{L(F)}$ and $Q(F)$ is defined in formula~\eqref{QF}.
\end{lem}
\begin{pf} 
Let $\A$ be the first component computed by Algorithm~\ref{MainRC}  and $p$ be the number of elements of $\A$. Note that
\begin{align}\label{eqHF}
p\Le n\cdot2^{H(\A)+m}.
\end{align}
We then
have
$$
\{F\} = [\A]:H_\A^\infty\cap\bigcap_{i=1}^p\{F,\i_i\}\cap\bigcap_{i=1}^p\{F,\s_i\}.
$$
If $1 \in (\A):H_\A^\infty$, we have 
$$
\s_1\cdot\ldots\cdot\s_p\cdot\i_1\cdot\ldots\cdot\i_p \in \sqrt{\left(\A\right)}.
$$
Therefore,
$$
f\cdot\s_1\cdot\ldots\cdot\s_p\cdot\i_1\cdot\ldots\cdot\i_p \in \sqrt{\left(\A\right)}.
$$
Consider now the case when $\A$ is as in Lemma~\ref{lemmaf} that gives us
$$
f\cdot\s_1^{j_1}\cdot\ldots\cdot\s_p^{j_p}\cdot\i_1^{k_1}\cdot\ldots\cdot\i_p^{k_p} \in \left(\A^{(\Le \ord f)}\right)
$$
for some non-negative integers $j_1,\ldots,j_p$
and $k_1,\ldots,k_p$.
Therefore,
$$
f\cdot\s_1\cdot\ldots\cdot\s_p\cdot\i_1\cdot\ldots\cdot\i_p \in \sqrt{\left(\A^{(\Le \ord f)}\right)}.
$$
Hence, by \cite[Corollary 1.7]{Kollar}, which gives an upper bound for degrees
in the algebraic Nullstellensatz when the degrees of the generating
polynomials are not equal to two, after squaring all elements of
$\A^{(\Le \ord f)}$ of degree two, we obtain that
\begin{align}\label{degreeinA}
(f\cdot\s_1\cdot\ldots\cdot\s_p\cdot\i_1\cdot\ldots\cdot\i_p)^d \in \left(\A^{(\Le \ord f)}\right),
\end{align}
where $d$ is defined by~\eqref{dformula}. Indeed,
 $n\cdot 2^{H(\A)+m+\ord f}$ bounds the number of algebraic indeterminates
in $\A^{(\Le\ord f)}$ and $f$, $H(\A)\Le H(F)\cdot2^{L(F)}$
 (because 
 at each step of
Algorithm~\ref{MainRC} the maximal order doubles at most), and by Corollary~\ref{maxcor} we have 
\begin{align*}
\max\left(4,D\left(\A^{(\Le \ord f)}\right)\right)& = \max(4,D(\A)) \Le\\
&\Le A\left(m+7,Q(F)-1\right).
\end{align*}
We also have
$$
1 \in \{F,\s_i\}:f\ \text{and}\ 1\in \{F,\i_i\}:f
$$
for all $i$, $1\Le i \Le p$. Let $G$ be $F,\i_i$ or $F,\s_i$  with the maximal $t(G,f)$. It then follows that
\begin{align*}
f^{k_j} = f_j + h_j,
\end{align*}
where $f_j \in \sqrt{\left(F^{\left(\Le t(G,f)\right)}\right)}$, $h_j \in \sqrt{\left(\i_i^{(\Le t(G,f))}\right)}$ or $\sqrt{\left(\s_i^{(\Le t(G,f))}\right)}$, and $k_j$
is a natural number.
Multiplying the above expressions, we obtain that
\begin{align}\label{Expressing1}
f^{\sum k_j} = g + h,
\end{align}
where $g \in \left(F^{\left(\Le t(G,f)\right)}\right)$ and
$$
h \in \left((\theta_1\s_1)\cdot\ldots\cdot(\theta_p\s_p)\cdot(\theta_1'\i_1)\cdot\ldots\cdot(\theta_p'\i_p) | \ord(\theta_k),\ord(\theta_l') \Le t(G,f)\right).
$$
Moreover, again since 
at each step of
Algorithm~\ref{MainRC} the maximal order doubles at most, we have
$$
\left(\A^{(\Le q)}\right)\subset \left({F^{(\Le q + H(F)\cdot 2^{L(F)})}}\right)
$$
for any $q \in \Z_{\Ge 0}$. By Lemma~\ref{degreelem} and inclusion~\eqref{degreeinA} we have
$$
(\theta_1\s_1)\cdot\ldots\cdot(\theta_p\s_p)\cdot(\theta_1'\i_1)\cdot\ldots\cdot(\theta_p'\i_p) \in \sqrt{\left(A^{\left(\Le \ord f + 4^{(2p+1)\cdot t(G,f)+1}d\right)}\right)}:f,
$$
where $d$ is defined in~\eqref{dformula}.
Hence, 
\begin{align}\label{InitialsInclusion}
(\theta_1\s_1)\cdot\ldots\cdot(\theta_p\s_p)&\cdot(\theta_1'\i_1)\cdot\ldots\cdot(\theta_p'\i_p) \in\notag\\
&\in \sqrt{\left(F^{\left(\Le \ord f + H(F)\cdot 2^{L(F)}+4^{(2p+1)\cdot t(G,f)+1}d\right)}\right)}:f
\end{align}
for all $\theta_k$ and $\theta'_k$ with $\ord(\theta_k),\:\ord(\theta'_k) \Le t(G,f)$, $1 \Le k \Le p$. Thus,
from inequality~\eqref{eqHF} and inclusion~\eqref{InitialsInclusion} it follows that
$$h \in \sqrt{\left(F^{\left(\Le \ord f + H(F)\cdot 2^{L(F)}+4^{\left(2n\cdot 2^{H(\A)}+1\right)\cdot t(G,f)+1}d\right)}\right)}:f,$$
which finishes the proof because we have representation~\eqref{Expressing1} and $g \in \left(F^{\left(\Le t(G,f)\right)}\right)$.
\end{pf}

\begin{thm}\label{MainTheor} We have
\begin{align}\label{NullBoundSolved}
t(F,f) \Le  A(m+8,\max(n,H(F\cup f),D(F\cup f))).
\end{align}
\end{thm}
\begin{pf} 
We begin the proof by recalling Corollary~\ref{maxcor}, which states that at
any stage of Algorithm~\ref{MainRC} the orders and degrees of differential
polynomials computed by this algorithm do not exceed the bound
$$E := A\left(m+7, Q(F\cup f)-1\right),$$
where the function $Q$ is defined in formula~\eqref{QF} and $Q(F\cup f) \Ge Q(F)$.
In particular,
$$\ord f\Le E,\;\; H(F)\Le E,$$
and by Proposition~\ref{MRCCorrect}
$$L(F)\Le \log_2 E.$$
It follows directly from the definition of $Q(F)$ 
that 
$$n\Le E,\;\; m\Le E,\; 100\Le E.$$
Using these inequalities, we can bound the quantity $d$ defined in \eqref{dformula} as
$$d\Le E^{E\cdot 2^{E^2+2E}}\Le E^{E^{E^E}},$$
whence 
$$d\cdot 4^d \Le 4^{2d} \Le E^{E^{E^{E^E}}}.$$
To simplify the latter formula, we note that 
$$\underbrace{k^{k^{k^{\ldots}}}}_{p{\rm\ times}} \Le 
\underbrace{2^{2^{2^{\ldots}}}}_{pk{\rm\ times}},$$
for all natural numbers $k$ and $p$, which can be easily derived by induction, using 
the inequality $ab\Le a^b$, which holds for all natural numbers $a,b \Ge 2$.
Thus we get:
$$d\cdot 4^d \Le 4^{2d} \Le \underbrace{2^{2^{2^{\ldots}}}}_{5E{\rm\ times}}.$$
This allows us to obtain the following inequality from \eqref{MainInequality}:
\begin{align*}
t(F,f)+1 &\Le E+1 + E^2 + \left(\underbrace{2^{2^{2^{\ldots}}}}_{5E{\rm\ times}}\right)^{t(G,f)+1} \Le \left(3\cdot \underbrace{2^{2^{2^{\ldots}}}}_{5E{\rm\ times}}\right)^{t(G,f)+1}\Le \\
&\Le
\left(\underbrace{2^{2^{2^{\ldots}}}}_{6E{\rm\ times}}\right)^{t(G,f)+1}.
\end{align*}
Now we use the fact that
$$A(4,k) = \underbrace{2^{2^{2^{\ldots}}}}_{k+3{\rm\ times}}-3$$
and simplify the above formula as
$$t(F,f)+1 \Le A(4,6E)^{t(G,f)+1},$$
and noting that $A(4,6E) \Le A(4,A(5,E-1))=A(5,E)$
yields 
$$t(F,f)+1 \Le A(5,E)^{t(G,f)+1}.$$

Now recall that the set $G$, defined in the proof of Lemma~\ref{NullBound}, is the set obtained from $F$ by Algorithm~\ref{MainRC} by adding to it an initial or separant. If now we take 
$$E_G := A\left(m+7, Q(G)-1\right),$$
and let $G_2$ be the set obtained from $G$ at the next iteration of the Algorithm~\ref{MainRC} by adding an initial or separant, we can similarly write
$$t(G,f)+1\Le A(5,E_G)^{t(G_2,f)+1}.$$
We continue recursively writing similar inequalities for $G_2,G_3,\ldots$, noting that the length of this chain does not exceed the number of iterations in Algorithm~\ref{MainRC}, that is, $L(F)\Le\log_2 E$. We note also that all quantities
$E, E_G, E_{G_2}, \ldots$ arising in these inequalities can be uniformly bounded by
$$\bar E := A\left(m+7, \max\left(9, n, 2^{9E},E\right)\right) = A\left(m+7,2^{9E}\right)\Le A(m+7,A(4,E)),$$
since the orders and degrees are uniformly bounded by $E$. Thus, we have
$$t(F,f)+1 \Le \underbrace{W^{W^{W^{\ldots}}}}_{E{\rm\ times}},$$
where $W = A(5,\bar E)$, which implies that
\begin{align*}
t(F,f)+1 &\Le \underbrace{2^{2^{2^{\ldots}}}}_{WE{\rm\ times}}\Le A(4,WE)\Le
A\left(4,W^2\right)\Le 
A(4,A(5,W-1))= \\
&= A(5,W)\Le A\left(5,A\left(6,\bar E-1\right)\right) = A\left(6,\bar E\right)=\\
&=A(6,A(m+7, A(4,E))) \Le A(m+6,A(m+7,A(4,E))) =\\
&=A(m+7,A(4,E)+1)) = \\
&=A(m+7,A(4,A(m+7,Q(F\cup f)-1))+1) \Le\\
&\Le A(m+7,A(m+6,A(m+7,Q(F\cup f)-1))+1)=\\
&=A(m+7,A(m+7,Q(F\cup f))+1). 
\end{align*}
Note that $$Q(F\cup f) \Le A(4,B) \Le A(m+8,B-1)-1,$$
where $B := \max(n,H(F\cup f),D(F\cup f)) \Ge 1$, since $n\Ge 1$.
Therefore, 
\begin{align*}
A(m+7,&A(m+7,Q(F\cup f))+1) \Le\\
&\Le A(m+7,A(m+7,A(m+8,B-1)-1)+1) \Le\\
&\Le A(m+7,A(m+7,A(m+8,B-1))-1) =\\
&=A(m+7,A(m+8,B)-1)=\\
&=A(m+8,B).
\end{align*}
\end{pf}

\section{Model-theoretic proof of existence of the bound}\label{MichaelsProof}
The following argument was shown to the authors by Michael Singer. In this
section we refer the reader for ultrafilters and construction of ultraproducts to books in model
theory, for instance \cite{Hodges,Dave}.
Let $\bar\K$ be the differential closure of $\K$ (see \cite[Definition 3.2]{PhyllisMichael} and the references given there) and $q \in \Z_{\Ge 0}$. 
We would like to emphasise
that in the statement below we had to {\it fix in advance} the number $r$ of differential
polynomials in $F$ to be able to use ultraproducts. So, in Theorem~\ref{modeltheorythm} the variable $r$ is quantified before the bounding function $\beta$. However, the constructive bound that we obtained in Theorem~\ref{MainTheor} {\it does not have} such a restriction, because it depends solely on the orders and degrees
in $F$ and $f$, but not on the number of elements in $F$.

\begin{thm}\label{modeltheorythm} For every $r \in \Z_{\Ge0}$ there exists a function $\beta : {\Z_{\Ge 0}}^3 \to \Z_{\Ge 0}$ such that for any
$q \in\Z_{\Ge 0}$ and $F \subset \K\{y_1,\ldots,y_n\}$ with 
$$
|F| = r,\ \ \max(H(F),D(F))\Le q,\ \ \text{and}\ \ 1\in [F]
$$
we have
$$
1 \in \left(F^{(\Le\beta(q,r,n))}\right).
$$ 
\end{thm}
\begin{pf}
Assume that the statement is wrong, that is, there exist $r,\:q\in\Z_{\Ge 0}$ such
that for any $\alpha \in \Z_{\Ge 0}$ there
exist $p_{1,\alpha},\ldots,p_{r,\alpha} \in \K\{y_1,\ldots,y_n\}$ with 
$$\max(H(p_{ij}),D(p_{ij}))\Le q$$ such that
$$1 \in [p_{1,\alpha},\ldots,p_{r,\alpha}]$$ 
and 
\begin{align}\label{onene}
1 \ne \sum\limits_{i=1}^r\sum\limits_{j=0}^\alpha q_{i,j}p_{i,\alpha}^{(j)}
\end{align}
for all $q_{i,j} \in \K\{y_1,\ldots,y_n\}$ of order less than or equal to $q+\alpha$. Again, it is essential
here that $r$ does not depend on $\alpha$.
For a maximal differential ideal $M$ in the differential ring $\prod_{i \in \Z_{\Ge 0}} \bar\K$
denote the differential ring $\left(\prod_{i \in \Z_{\Ge 0}} \bar\K\right)/M$ by $K_M$.
There is a natural
differential ring homomorphism
$$
 \left(\prod_{i \in \Z_{\Ge 0}} \bar\K\right)\{y_1,\ldots,y_n\} \to K_M\{y_1,\ldots,y_n\} =: R.
$$
We shall now make a special choice of the maximal differential ideal $M$. Let
$\mathcal{F}$ be the filter consisting of all cofinite subsets of $\Z_{\Ge 0}$.
Then, there exists an ultrafilter $\mathcal{U}$ containing $\mathcal{F}$. Since
the field $\bar\K$ is differentially closed, by
\L o\'s' theorem \cite[Theorem 8.5.3]{Hodges} the ultraproduct 
$$K := \prod \bar\K/\mathcal{U}$$ is a differentially closed
field with the following property. 

Let $\bar{a} = (a_0,a_1,a_2,\ldots)$ and $\bar{b} = (b_0,b_1,b_2,\ldots) \in K$. Then, we have: if $\bar{a} = \bar{b}$ then $a_i = b_i$
for infinitely many indices $i$. We now take $M$ to be the kernel of the differential
ring homomorphism
$$
\prod_{i\in\Z_{\Ge 0}}\bar\K \to K.
$$
Let $\bar{p_i}$ be the image of $(p_{i,1},p_{i,2},p_{i,3},\ldots)$ in $R$. This is defined correctly as all $p_{i,j}$ have orders and degrees bounded. 
Assume that $(z_1,\ldots,z_n) \in (K_M)^n$ is a zero of $\bar{p_i}$ for all $i$. Then,
for each $i$, $1\Le i\Le r$, there exists $V_i \subset \Z_{\Ge 0}$, $V_i \in \mathcal{U}$,
such that 
$$p_{i,j}(z_{1,j},\ldots,z_{n,j}) = 0$$ 
for all $j \in V_i$,
where $(z_{t,1},z_{t,2},z_{t,3},\ldots)$ is mapped to $z_t$ under the mentioned differential ring homomorphism for each $t$, $1\Le t \Le n$. Since $V_1\cap\ldots\cap V_r \in \mathcal{U}$ and $\varnothing \notin \mathcal{U}$, there is an index 
$$
j \in V_1\cap\ldots\cap V_r.
$$ Therefore,
$$
p_{1,j}(z_{1,j},\ldots,z_{n,j}) = \ldots = p_{r,j}(z_{1,j},\ldots,z_{n,j}) = 0.
$$
Since $\bar\K$ is differentially closed, this contradicts to $1 \in [p_{1,j}, \ldots, p_{r,j}]$ in the differential ring $\bar\K\{y_1,\ldots,y_n\}$.
Therefore, since the field $K_M$ is differentially closed, we have $1 \in [\bar{p_1},\ldots,\bar{p_r}]$. Hence, there exist $\gamma \in \Z_{\Ge 0}$ and differential polynomials
$\bar{q_{ij}} \in K_M\{y_1,\ldots,y_n\}$ with $\ord\bar{q_{ij}} < \gamma + q$ so that
$$
1 = \sum_{i=1}^r\sum_{j=0}^\gamma\bar{q_{ij}}\bar{p_i}^{(j)}.
$$
Again, due to our choice of $M$ (that is, due to the fact
that $\mathcal{U}$ is an ultrafilter), there exists $\alpha \in \Z_{\Ge 0}$ with $\alpha > \gamma$ such that 
$$
1 = \sum_{i=1}^r\sum_{j=0}^\gamma q_{ij}p_{i,\alpha}^{(j)},
$$
where 
$q_{i,j} \in \bar\K\{y_1,\ldots,y_n\}$ of order less than $\alpha+q$. Since $p_{i,\alpha} \in \K\{y_1,\ldots,y_n\}$ for all $i$, $1\Le i\Le r$, by taking
a basis of $\bar\K$ over $\K$ we may assume that in fact $q_{i,j}\in \K\{y_1,\ldots,y_n\}$ for all $i$ and $j$, $1\Le i \Le r$, $0\Le j\Le \gamma$. This contradicts to~\eqref{onene}. Thus, our
initial assumption was wrong.
\end{pf}

\section{Conclusions} We have obtained the first bound on orders for
the differential Nullstellensatz. Surely, one can improve the bound and find many applications of it.
A general programme which is being realized here is as follows. The differential elimination algorithms
would be very useful for applications if there
were faster versions of them. Our work on bounding
orders could lead to:
\begin{enumerate}
\item understanding complexity estimates for the differential elimination,
\item developing combined and separated differential and high performance algebraic algorithms.
\end{enumerate}
One of the ideas is, instead of using the usual
differential elimination, perform all differentiations at the beginning of the process and then use only fast algebraic methods. We hope that our bounds will contribute to this programme.

\section{Acknowledgements} We thank Michael Singer for very helpful
comments, support, and for the model theoretic proof of the
differential Nullstellensatz. We are grateful to Daniel Bertrand for
encouraging us to solve the problem. We appreciate the help of Erich
Kaltofen, Teresa Krick, Alice Medvedev, and Eric Schost in finding references to the previous
work on the algebraic version of the effective Nullstellensatz,
on bounds for the lengths of monomial sequences, and on model theory. 
We are  grateful to the referees for important suggestions.

\bibliographystyle{elsart-num}
\bibliography{diffnull}
\end{document}